\documentclass[a4paper,12pt]{article}
\usepackage[dvips]{epsfig}
\usepackage{amsmath, amssymb, amsbsy, amstext, amscd,color,enumerate,amsfonts}
\usepackage{euscript, graphicx, hyperref, srcltx, tikz, indentfirst}

\input amssym.def
\input amssym.tex

\usepackage{listings}

\lstdefinelanguage{Magma}{
  keywords={for, while, if, then, else, function, return}, 
  sensitive=true,
  morecomment=[l]{//},
  morecomment=[s]{/*}{*/},
  morestring=[b]',
  morestring=[b]"
}

\newtheorem{lem}{Lemma}[section]%
\newtheorem{theorem}[lem]{Theorem}%
\newtheorem{cor}[lem]{Corollary}%
\newtheorem{que}[lem]{Question}%
\newtheorem{prop}[lem]{Proposition}%

\def\og{\overline G}  \def\ox{\overline X} 
\def\di{\bigm|} \def\lg{\langle} \def\rg{\rangle} \def\ola{\overline a} \def\olb{\overline b}
\def\PSL{\hbox{\rm PSL}} \def\Aut{\hbox{\rm Aut}~} \def\PGL{\hbox{\rm PGL}} \def\O{\Omega}
\def\AGL{\hbox{\rm AGL}} \def\GL{\hbox{\rm GL}} \def\character{ ~ {\rm char} ~} \def\o{{\rm o}}
\def\ZZ{\mathbb{Z}} \def\f{\noindent} \def\qed{\hfill $\Box$} \def\demo{\f {\bf Proof}\hskip10pt}

\begin{document}
\begin{center}
{\bf\large  The product of nonabelian simple groups and dihedral groups}
\end{center}

\begin{center}
Hao Yu{\small \footnotemark}\\
\medskip
{\it {\small
School of Mathematical Sciences, Capital Normal University,\\
Bejing 100048, People's Republic of China}}
\end{center}

\footnotetext{Corresponding author: haoyu@gxu.edu.cn.
This work was supported in part  by the National Natural Science Foundation of China
 (12471332).}

\begin{abstract}
  Let $X=GD$ be a group, where $G$ is a nonabelian simple group and $D$ is a dihedral group.
  These groups $X$ are closely related to regular Cayley maps.
  The main theorems of this paper describes $X$.
\end{abstract}

\renewcommand{\thefootnote}{\empty}

\section{Introduction}
Throughout the paper all groups are finite.
A group $X$ is said to be properly {\it factorizable} if $X=GH$ for two proper subgroups $G$ and $H$ of $X$, while the expression $X=GH$ is called a {\it factorization} of $X$.
Also, $X=GH$ is called a {\it product group} of $G$ and $H$.
Furthermore, if $G\cap H=1$, then  we say that $X$ has an {\it exact factorization}.

Factorizations of groups naturally arise from the well-known Frattini's argument, including its version in permutation groups.
One of the most famous results about factorized groups might be one of theorems of It\^o, saying that any group is metabelian whenever it is the product of two abelian subgroups (see \cite{Ito1}).
Later, Wielandt and  Kegel showed that the product of two nilpotent subgroups must be soluble (see \cite{Wie1958} and \cite{Ke1961}).
Douglas showed that  the product of two cyclic groups must be super-solvable (see \cite{D1961}).
In \cite{LPS1990} the maximal factorizations of finite almost simple groups were determined  and
the  factorizations of almost simple groups with a solvable factor were determined in \cite{LX2022}.
There are many other papers related to factorizations, for instance, finite products of soluble groups, factorizations with one nilpotent factor and so on.
Here we are not able to list all references and the readers may refer to a survey paper \cite{AK2009}.

For a factorization $X=GH$, let $\O=[X : G]$ be the set of right cosets of $G$ in $X$.
Suppose that $G$ is core free. Then $X$ acts faithfully on $\O$.
Hence, $X$ is a permutation group with degree $|\O|$ containing a transitive subgroup $H$.
The study of a permutation group containing a transitive subgroup played an important role in the history of permutation group theory.
A transitive permutation group $X$ is called {\it primitive} if $\O$ has no nontrivial $X$-invariant partitions.
Burnside and Schur proved that the primitive permutation group containing a cyclic regular subgroup of order $p^m$,
where $m>1$, or composite order is 2-transitive (see \cite{B1911,S1933}).
In 1949,  Wielandt showed that such group containing a dihedral regular subgroup is 2-transitive \cite{W1947},
and in 1957, Scott proved that such group containing a regular generalized quaternion subgroup is 2-transitive \cite{S1957}.
A permutation group is called {\it quasiprimitive} if each of its non-trivial normal subgroups is transitive.
Li classified quasiprimitive permutation groups containing a dihedral regular subgroup \cite{L2006}, and recently, Li, Pan and Xia characterized  quasiprimitive permutation groups containing a metacyclic transitive subgroup \cite{LPX2021}.
In this paper, we shall determine permutation groups that contain a dihedral transitive subgroup
and a point stabilizer isomorphic to a nonabelian simple group.

This paper is also motivated from Cayley maps, see \cite{J1997}.
The automorphism group of an orientably (unoriented) regular Cayley map of $G$ is a product group of $G$ and a cyclic (dihedral) group, separately.
Thus, classifying orientably (respectively, unoriented) regular Cayley maps of $G$ is related to determining product groups of $G$ and a cyclic (respectively, dihedral) group.
Since orientably regular Cayley maps of nonabelian simple groups are characterized in \cite{CDL2022}, we shall consider their unoriented regular Cayley maps.
Hence, we need to study product groups of  nonabelian simple groups and dihedral groups.

The following theorem determines the product of nonabelian simple groups and dihedral groups.
Remind that a group $G=K.H$ denotes an extension of a group $K$ by a group $H$.
the core $G_X$ of a subgroup $G$ in a group $X$ is $G_X=\cap_{x\in X}G^x$, and when we say that $G$ is core free in $X$, it means $G_X=1.$

\begin{theorem} \label{main1}
Let $X=GD$ be a product group of a nonabelian simple group $G$ and a dihedral group $D$.
Then either $G \character X$ or $G_X=1$ and one of the following holds:
\begin{enumerate}
  \item[\rm(1)] $G\cap D=1$, and either
  \begin{enumerate}
  \item[\rm(1.1)] $D_X=1$ and the triple $(X,G,D)$ is given in Table~\ref{dihedral-skew}; or
  \item[\rm(1.2)] $D_X\neq1$ is cyclic and $X=D_X\rtimes Y$, where $(Y,G)$ is
  either $(A_{2m+1}\rtimes\ZZ_2,A_{2m})$, where $m\ge 3$, or $(PGL(2,11),A_5)$.
  \end{enumerate}
  \item[\rm(2)] $G\cap D\cong\ZZ_2$, and $(X,G,D)$ is either $(A_{m+1}, A_m,D_{2(m+1)})$ or $(A_{m+1}\times\ZZ_2, A_m,D_{4(m+1)})$, with $4\di m$ and $m\ge 8$.
\end{enumerate}
\end{theorem}

\begin{table}\centering  \caption{The product group $X$ of $G$ and $D$}\label{dihedral-skew}
  \begin{tabular}{llll}\hline
   Row &      $X$     & $G$ &  $D$\\
  \hline
   1   & $\AGL(3,2)$  &$GL(3,2)$   &$D_8$\\
   2   & $M_{12}$     & $M_{11}$  &$D_{12}$\\
   3   & $M_{24}$     & $M_{23}$      &$D_{24}$\\
   4   & $A_{4m}, m\ge 2$     &$A_{4m-1}$  &$D_{4m}$\\
   5   &$PGL(2,11)$   & $A_5$     &$D_{22}$ \\
   6   &$A_{2m+1}\rtimes\ZZ_2, m\ge 3$ & $ A_{2m}$    &$D_{2(2m+1)}$\\
   7   &$\Aut(M_{12})$& $M_{11}$       &$D_{24}$\\
   8   &$A_{4m}\rtimes\ZZ_2,  m\ge 2$ & $A_{4m-1}$  &$D_{8m}$\\
  \hline
  \end{tabular}
\end{table}

\vskip 3mm
\f{\it Remark:} $X$ acts quasiprimitively on $[X:G]$  in Rows 1-4 of Table~\ref{dihedral-skew} and $X$ is not quasiprimitive in Rows 5-8.

\section{Preliminaries}
Notations and terms used in the paper are standard, see  \cite{R1995}.
In what follows, we introduce some known group theoretical results for later
use.

\begin{prop}\cite[Theorem 1.5 and Lemma 7.2]{L2006}\label{Lprimitive}
Let $X$ be a quasiprimitive permutation group acting on $\O$ of degree $n$.
Then, $X$ contains a dihedral regular subgroup $G$ if and only if $X$ is 2-transitive, and one of the following holds, where $\omega \in \Omega$:
\begin{enumerate}
  \item[\rm(i)] $(X, G, X_\omega)=(A_4, D_4,\ZZ_3), (S_4, D_4,D_6), (\AGL(3,2), D_8,\GL(3,2)),
  (\AGL(4,2), D_{16},$ $\GL(4,2)), (\ZZ_2^4\rtimes A_6, D_{16},A_6)$, or $(\ZZ_2^4\rtimes A_7, D_{16},A_7)$;
  \item[\rm(ii)] $(X, G, X_\omega)=(M_{12}, D_{12}, M_{11}),(M_{22} .2, D_{22}, \operatorname{PSL}(3,4) .2)$, or $(M_{24}, D_{24}, M_{23})$;
  \item[\rm(iii)] $(X, G, X_\omega)=(S_{2 m}, D_{2 m}, S_{2 m-1})$ or $(A_{4 m}, D_{4 m}, A_{4 m-1})$;
  \item[\rm(iv)] $X=\PSL(2, p^e).O, G=D_{p^e+1}$, and
  $X_\omega \unrhd \ZZ_p^e \rtimes \ZZ_{\frac{p^e-1}{2}}.O$,
  where $p^e \equiv3\pmod 4$, and $O \leq {\rm Out}(\PSL(2, p^e)) \cong \ZZ_2 \times \ZZ_e ;$
  \item[\rm(v)] $X=\PGL(2, p^e) \ZZ_f, G=D_{p^e+1}$, and $X_\omega=\ZZ_p^e \rtimes \ZZ_{p^e-1}$, where $p^e \equiv 1$ $\pmod 4$, and $f \mid e$.
\end{enumerate}
\end{prop}

\vskip 3mm
Picking  up all point stabilizers which are nonabelian simple groups in  Proposition~\ref{Lprimitive} and using Magma \cite{BCP1997}, we immediately get the following corollary.
\begin{cor}\label{primitive}
Let $X=GD$ be a product group of a nonabelian simple group $G$ and a dihedral group $D$ such that $G\cap D=1$. Suppose that $X$ acts faithfully and quasiprimitively on $[X:G]$.
Then $(X, G, D)=(\AGL(3,2),GL(3,2),D_8)$, $(M_{12},M_{11},D_{12})$, $(M_{24}, M_{23}, D_{24})$,
or  $(A_{4m},A_{4m-1},D_{4m})$ with $m\ge2$.
\end{cor}
\demo
By hypothesis $X$ is quasiprimitive permutation group with a dihedral regular subgroup $D$ and its stabilizer is isomorphic to the nonabelian simple group $G$.
By Proposition~\ref{Lprimitive}, we have:
$(X, D, G)=(\AGL(3,2), D_8, \GL(3,2)), (\AGL(4,2), D_{16},\GL(4,2)), (\ZZ_2^4\rtimes A_6, D_{16}, A_6)$,
$(\ZZ_2^4\rtimes A_7, D_{16}, A_7)$, $(M_{12}, D_{12}, M_{11})$, $(M_{24}, D_{24}, M_{23})$ or
$(A_{4 m}, D_{4 m}, A_{4 m-1})$, where $m\ge2$.
Reminding $D\cap G=1$ and using Magma \cite{BCP1997}, we get $X\notin\{ \AGL(4,2), \ZZ_2^4\rtimes A_6$, $\ZZ_2^4\rtimes A_7\}$.
Therefore, we get that $(X, G, D)=(\AGL(3,2),GL(3,2),D_8)$, $(M_{12},M_{11},D_{12})$, $(M_{24}, M_{23}, D_{24})$ or  $(A_{4m},A_{4m-1},D_{4m})$ with $m\ge2$, as desired.
\qed

\begin{prop}\cite[Corollary 1.2]{BCV2022}\label{skew}
Let $X=GC$ be a group where $G$ is a nonabelian simple group and $C$ is cyclic
such that $G\cap C=1$. Suppose that $C$ is core-free.
Then either 
\begin{enumerate}
  \item[\rm(1)] $X=G\rtimes C$ and $C\leq\Aut(G)$; or
  \item[\rm(2)] $(X, G)=(\PSL(2,11), A_5),(M_{23}, M_{22})$, or $(A_{m+1}, A_m)$ with $m \ge 6$ even.
\end{enumerate}
\end{prop}

\vskip 3mm
For a group $H$, consider $G:=K.H$.
In particular, $G$ is a {\it proper central extension} of $K$ by $H$ if $K \leq Z(G)\cap G'$ \cite{S1982}.
Then a central extension $G$ is a proper central extension, $K\times H$ or a mix of both.

\begin{prop}\cite[Proposition 2.3]{CDL2022}\label{Schur}
For simple groups $\operatorname{PSL}(2,11), M_{23}$ and $A_n$ for $5 \leqslant n \neq 6$, we have the following results:
\begin{enumerate}
  \item[\rm(1)] the proper central extensions of $\operatorname{PSL}(2,11)$ (or $A_n$) don't contain $A_5$ (or $A_{n-1}$); and
  \item[\rm(2)] there exists no nontrivial proper central extension of $M_{23}$.
\end{enumerate}
\end{prop}

\vskip 3mm
Using Magma \cite{BCP1997}, we get the following facts.
\begin{lem}\label{Magma}
\begin{enumerate}
  \item[\rm(1)] If $\PSL(2,11)\rtimes\ZZ_2$ contains a subgroup $D_{22}$, then $\PSL(2,11)\rtimes\ZZ_2=\PGL(2,11)=A_5D_{22}$;
  \item[\rm(2)] if $M_{12}\rtimes\ZZ_2$ contains a subgroup $D_{24}$, then $M_{12}\rtimes\ZZ_2=\Aut(M_{12})=M_{11}D_{24}$;
  \item[\rm(3)] there exists no $\AGL(3,2)\rtimes\ZZ_2$ containing a subgroup $D_{16}$; and
  \item[\rm(4)] there exists no $M_{23}\rtimes\ZZ_2$ (resp. $M_{24}\rtimes\ZZ_2$) containing a subgroup $D_{46}$ (resp. $D_{48}$).
\end{enumerate}
\end{lem}

\section{Proof of Theorem~\ref{main1}}
Throughout this section, let $X=GD$ be a product group of a nonabelian simple group $G$ and a dihedral group $D$.
Set $D=\lg a\rg \rtimes\lg b\rg$, where $\langle b \rangle$ has order 2.
For a convenience,  {\it the product group of a nonabelian simple group and a dihedral group} is  abbreviated as {\it PSDG}.
In other words, a group $X$ is {\it PSDG} if $X$ has a factorization $X=GD$ with a nonabelian simple group $G$ and a dihedral group $D$.

\vskip 3mm 
Since $G_X\unlhd G$ and $G$ is a nonabelian simple group, we get that $G_X$ is either $G$ or 1.
If $G_X=G$, that is $G\unlhd X$, then $G\character X$ as $G$ is a nonabelian  simple group and $G/G_X$ is soluble, which is the former case  of Theorem~\ref{main1}.
So in what follows, we assume $G_X=1$.

\begin{lem}\label{D-goup1}
Suppose $G \le J \le X$. Then $J = G(J \cap D)$ and $J_X = \bigcap_{d \in D} J^d.$
In particular, if $E \le D \cap J$ is normal in $D$, then $E \le J_X$.
\end{lem}

\demo
Let $G \le J \le X$.
Then by the Dedekind modular law, $J = G(J \cap D)$, which implies $J^{g}$ for any $g\in G$.
Thus,
$$J_X = \bigcap_{x \in X} J^x= \bigcap_{g \in G,d\in D} J^{gd} = \bigcap_{d \in D} J^d,$$
as desired.

Let $E \le D \cap J$ be normal in $D$.
Then 
$$E=\bigcap_{d \in D} E^d\leq \bigcap_{d \in D} J^d=J_X,$$
as desired.
\qed

\begin{lem}\label{D-goup}
Suppose that $G$ is core-free.
Then $G\cap D$ is either 1 or $\lg b_1\rg$ where $b_1\in D\setminus\lg a\rg$.
\end{lem}
\demo
Since $G_X=1$, $D \cap G$ contains no non-trivial normal subgroup of $D$ by Lemma~\ref{D-goup1}. 
Hence $G \cap D$ has order 1 or 2.
In particular, if $G \cap D$ has order 2, then $G\cap D=\lg b_1\rg$ with $b_1\in D\setminus\lg a\rg$, as desired.
\qed

\vskip 3mm
\f{\it Remark:}
In fact, if we remove the condition that $G$ is a nonabelian group, then the above lemma also holds.

\vskip 3mm
By Lemma~\ref{D-goup}, we get $|G\cap D|\leq2$ for $X=GD$.
Depending on $|G\cap D|$, we divide the proof into the following two subsections.

\subsection{ $G \cap D = 1$ }\label{subsec:GD1}

In this subsection, let $X=GD$ be {\it PSDG} such that $G\cap D=1$ and $G_X=1$.
Consider $D_X$, the core of $D$ in $X$.
Firstly, we deal with the case that $D_X=1$.
Before giving a classification,  the following useful lemma is needed.

\begin{lem}\label{Ga}
Suppose that $D_X=1$ and there exists a maximal subgroup $H$ of $X$ such that $G<H$ and $H\cap D\leq\lg a\rg$.
Then $X=(G\lg a\rg)\rtimes\lg b\rg$ and $(X,G,D)$ is either $(\PGL(2,11), A_5,D_{22})$ or $(A_{m+1}\rtimes\ZZ_2, A_{m},D_{2(m+1)})$ with $m\ge6$ even.
Moreover, $G\lg a\rg\character X$.
\end{lem}
\demo
Set $H \cap D = \langle a_1 \rangle \le \langle a \rangle$. 
Then  $\langle a_1 \rangle \le H_X$ by Lemma~\ref{D-goup1}.
Then we get $H_X=(H_X\cap G)\lg a_1\rg$.
Noting that $G$ is a nonabelian simple group and $H_X\cap G\unlhd G$,
we know that $H_X\cap G$ is either 1 or $G$.
For a contradiction, assume that $H_X\cap G=1$.
Then one can get $H_X=\lg a_1\rg\leq D_X=1$, a contradiction.
Therefore, $H_X\cap G=G$ and so $H=H_X\lhd X$.
By the maximality of $H$, we get $\lg a_1\rg=\lg a\rg$, which implies
$$H=G\lg a\rg\quad{\rm and}\quad X=H\rtimes\lg b\rg.$$
Moreover, we get
$$\lg a\rg_H=\cap_{a^i g\in H}\lg a\rg^{a^ig}=\cap_{b^ja^ig\in X}\lg a\rg^{b^ja^ig}=\lg a\rg_X\leq D_X=1.$$
Noting that $G\cap\lg a\rg=1$, by Proposition~\ref{skew}, we get that either $H=G\rtimes\lg a\rg$, or
$(H, G)=(\PSL(2,11), A_5)$, $(M_{23}, M_{22})$ or $(A_{m+1}, A_m)$ with $m\ge6$ even.
If $H = G \rtimes \langle a \rangle$, then $G$ is characteristic in $H$ and so normal in $X$, a contradiction.
Hence, $$(H, G,\lg a\rg)=(\PSL(2,11), A_5,\ZZ_{11}),(M_{23}, M_{22},\ZZ_{23})\,~{\rm or}~\,(A_{m+1}, A_m,\ZZ_{m+1}).$$
Since $X = H \rtimes \langle b \rangle \cong H \rtimes Z_2$ we now obtain  
$$(X, G,D)\in\{(\PSL(2,11)\rtimes\ZZ_2, A_5,D_{22}),(M_{23}\rtimes\ZZ_2, M_{22},D_{46}),(A_{m+1}\rtimes\ZZ_2, A_m,D_{2(m+1)})\}.$$

By Lemma~\ref{Magma} the first case gives $(X,G,D) = (PGL_2(11), A_5, D_{22})$ while the second case gives no examples. In the third family, we have $(X,G,D) = (A_{m+1} \rtimes Z_2, A_m, D_{2(m+1)})$ with $m\ge6$ even and we can realize this with $A_{m+1} \rtimes Z_2$ as follows.
    
Let $$D=\lg a,b\di a=(1\, 2 \cdots m+1),~b=(1\,m+1)(2\,m)\cdots(\frac{m}2\,\frac{m}2+2)(m+2\,m+3)\rg.$$
Then we get $D\cong D_{2(m+1)}$, $a\in A_{m+1}$ and $D\leq S_{m+3}$.
Let $G=A_m\leq S_m$ and $X=A_{m+1}\rtimes\lg b\rg$.
Then one yields $X=GD$.
Since $A_{m+1}$ is simple for $m\ge4$, we know $\lg a\rg_X=1$, which implies $D_X=1$.
Hence, $(X,G,D)$ can be $(A_{m+1}\rtimes\ZZ_2, A_m,D_{2(m+1)})$ with $m\ge6$ even, as desired.
Moreover, 
From the above construction, it is easy to see that $X\cong A_{m+1}\times\ZZ_2$ if $\frac m2$ is even, and $X\cong S_{m+1}$ if $\frac m2$ is odd.

Now, $(X,G,D)$ is either $(\PGL(2,11), A_5,D_{22})$ or $(A_{m+1}\rtimes\ZZ_2, A_m,D_{2(m+1)})$ with $m\ge6$ even.
Then we get that $(X,G\lg a\rg)$ is either $(\PGL(2,11), \PSL(2,11))$ or $(A_{m+1}\rtimes\ZZ_2, A_{m+1})$ with $m\ge6$ even.
Hence, $G\lg a\rg\character X$, as desired.
\qed

\begin{lem}\label{GX1}
Suppose that $D_X = 1$. Then $(X, G, D)$ is given in Table~\ref{dihedral-skew}.
\end{lem}
\demo
Since $G_X = 1$, the group $X$ acts faithfully on $\O=[X:G]$.
Suppose that $|D|=4$.
Then $X\leq S_4$, but there exists no nonabelian simple group in $S_4$, a contradiction.
So in what follows, we assume that $|D|>4$.
If $X$ is quasiprimitive, then by Corollary~\ref{primitive}, we get
$(X, G, D)=(\AGL(3,2),GL(3,2),D_8)$, $(M_{12},M_{11},D_{12})$,
$(M_{24}, M_{23}, D_{24})$ or $(A_{4m},A_{4m-1},D_{4m})$ with $m\ge2$.
This gives Rows 1--4 of Table~\ref{dihedral-skew}.
So in what follows, we may assume that $X$ is not quasiprimitive.

Since $X$ is not quasiprimitive, there exists a non-transitive normal subgroup $N$ of $X$. 
Since $N$ is not transitive $GN < X$.
Let $H$ be a maximal subgroup of $X$ containing $GN$.
Clearly, we know $N\leq H_X\neq1$.
By $G<H$, we get $H=G(H\cap D)$ where $H\cap D\neq1$.
If $H\cap D\leq\lg a\rg$, then by Lemma~\ref{Ga}, we get
$$(X,G,D)=(\PGL(2,11), A_5,D_{22}),(A_{m+1}\rtimes\ZZ_2, A_{m},D_{2(m+1)}),
~{\rm with}~m\ge6~{\rm even}.$$
This gives Rows 5 and 6 of Table~\ref{dihedral-skew}.
So in what follows, we assume that $H\cap D\not\leq\lg a\rg$.
In fact, we can make a more precise assumption  $\mathfrak{A}$:
for any maximal subgroup $L$ of $X$, if $G< L$, then $L\cap D\not\leq\lg a\rg.$

Set $H\cap D=\lg a_1,b_1\rg$ for some $\lg a_1\rg<\lg a\rg$ and $b_1\in D\setminus\lg a\rg$.
Note that $D=\lg a,b\rg=\lg a,b_1\rg=\lg a\rg\rtimes\lg b_1\rg$.
We may as well choose notations so that $b_1 = b$. Then $H=G\lg a_1,b\rg$.
By Lemma~\ref{D-goup1},  $\lg a_1\rg\leq H_X$ and $H_X\cap G$ is either 1 or $G$.
Depending on $H_X\cap G$, we can divide it into the following two cases.

\vskip 3mm
{\it Case 1: $H_X\cap G=1$.}
\vskip 3mm

Suppose that  $b\in H_X$. 
Then $\lg a_1,b\rg\leq H_X$.
Since $H_X\cap G=1$, we know $H_X=\lg a_1,b\rg\leq D_X=1$, a contradiction.
Hence, $b\notin H_X.$

In $\overline{H}=H/H_X=\overline{G\lg a_1,b\rg}=\overline{G\lg b\rg}$, 
we get $\overline{G}\cong G/(G\cap H_X)\cong G$ and $\lg \overline b\rg\cong \ZZ_2$.
If $\overline{G}\cap \lg \overline b\rg=1$, then $|H_X|=|\lg a_1\rg|$, which implies $H_X=\lg a_1\rg$,
and so $a_1=1$.
Then we assume that $\lg \overline b\rg\leq \overline{G}$.
Then $|H_X|=2|\lg a_1\rg|$ and so $H_X\cong \ZZ_{\o(a_1)}.\ZZ_2$.
Set $\mho_1(H_X)=\lg h^2\di h\in H_X\rg$.
Then by $(gb)^2\in\lg a_1\rg$, we get $\lg a_1^2\rg\leq\mho_1(H_X)\leq\lg a_1\rg$ and $\mho_1(H_X)\character H_X$.
Hence, $\lg a_1^2\rg\character H_X\lhd X$, which implies $\lg a_1^2\rg\lhd X$,
and so $\lg a_1^2\rg\leq D_X=1$.
We claim that $a_1=1$.
Indeed, otherwise, $H_X\cong D_4$ or $\ZZ_4$ as $a_1$ is an involution.
Since $H/C_H(H_X)\leq \Aut(H_X)$, $H=GH_X$, $H_X\leq C_H(H_X)$ and 
$GC_H(H_X)/C_H(H_X)$ is solvable, we get $C_H(H_X)=H$, which implies $\lg G,a\rg\leq C_X(\lg a_1\rg)$.
Then $\lg a_1\rg\lhd X$ and so $\lg a_1\rg\leq D_X=1$, a contradiction.

Therefore, we get $$H_X=\lg gb\rg\cong\ZZ_2\,\,{\rm and}\,\,H=G\rtimes\lg b\rg=G\times\lg gb\rg,$$
where $g\in G.$
We claim $H\ntrianglelefteq X$.
If not, we get $G\character H\lhd X$, which implies $G\lhd X$, a contradiction.
Consider the quotient group $\ox:=X/H_X=\og\lg \ola\rg$.
Then $\og\ntrianglelefteq\ox$, which implies $\og_{\ox}=1$.
Set $\lg \ola_2\rg=\lg \ola\rg_{X/H_X}$.
Note that $H\leq G\lg a_2,b\rg<X$ and $\og\times \lg \ola_2\rg\leq\ox$.
By the maximality of $H$, we get $a_2=1$, that is $\lg \ola\rg_{X/H_X}=1$.
Hence, $\lg\ola\rg$ is core-free in $X/H_X$.
By Proposition~\ref{skew}, we know
$(\ox, \og)=(\PSL(2,11), A_5),(M_{23}, M_{22})~{\rm or}~(A_{m+1}, A_{m}),$
where $m\ge6$ is even.
Hence, $(X, G)\in\{(\ZZ_2.\PSL(2,11), A_5),(\ZZ_2.M_{23}, M_{22}),(\ZZ_2.A_{m+1}, A_{m})\}.$
By Proposition~\ref{Schur} and Lemma~\ref{Magma}, we know that if $\ZZ_2.A_{m+1}= A_{m} D_{2(m+1)}$, then $\ZZ_2.A_{m+1}= A_{m+1}\times\ZZ_2$, and there exists no $\ZZ_2.\PSL(2,11)$ ($\ZZ_2.M_{23}$) such that $\ZZ_2.\PSL(2,11)=A_5 D_{22}$ ($\ZZ_2.M_{23}=M_{22} D_{46}$), separately.
Hence, we get $(X, G, D)=(A_{m+1}\times \ZZ_2, A_{m}, D_{2(m+1)}).$
But $G\lg a\rg=A_{m+1}$ is a maximal subgroup of $X$ and $G\lg a\rg\character X$, contradicting to the assumption $\mathfrak{A}$.

\vskip 3mm
{\it Case 2: $H_X\cap G=G$.}
\vskip 3mm

Clearly, we have $H_X=H_X\cap H=G(H_X\cap\lg a_1,b\rg)$.
Noting that $\lg a_1\rg\leq H_X$, we get that $H_X$ is either $G\lg a_1\rg$ or $G\lg a_1,b\rg$.
But, if $H_X=G\lg a_1\rg$, then $H_X\lg a\rg=G\lg a\rg$ is a maximal subgroup of $X$ and $G\lg a\rg\lhd X$,
contradicting to $\mathfrak{A}$.
Hence, $H_X=G\lg a_1,b\rg=H$, which implies $H\lhd X$.
Since $b\in H\unlhd X$ and $b^a=a^{-2}b$, we get $a^2\in H$.
Then $H=G\lg a^2,b\rg$, and so $X=H\rtimes\lg ab\rg$,
which implies $4\di|D|$.
Hence, $|D|\ge8$ as $|D|>4$.
Since $a^2\neq1$, we know that $H=G\lg a^2,b\rg$ is {\it PSDG} and $G\cap\lg a^2,b\rg=1$.
We claim that $G\ntriangleleft H$.
Otherwise, since $G$ is a nonabelian simple group, we get $G\character H\lhd X$, which implies $G\lhd X$, contradicting to $G_X=1$.
Therefore, $G_H=1$, which implies that $H$ acts faithfully on $[H:G]$.

Suppose that $H$ is quasiprimitive.
By Corollary~\ref{primitive}, we get
$(H, G, \lg a^2,b\rg)=(\AGL(3,2),\\ GL(3,2),D_8),\,(M_{12},M_{11},D_{12}),\,(M_{24}, M_{23}, D_{24}),\,(A_{4m},A_{4m-1},D_{4m}),$
with $m\ge2$, which implies
$(X, G, D)=(\AGL(3,2)\rtimes\ZZ_2,GL(3,2),D_{16}),\,(M_{12}\rtimes\ZZ_2,M_{11},D_{24})$, $(M_{24}\rtimes\ZZ_2, M_{23},\\ D_{48})$ or $(A_{4m}\rtimes\ZZ_2,A_{4m-1},D_{8m}).$
By Lemma~\ref{Magma}, we get
$$(X, G, D)\in\{(\Aut(M_{12}),M_{11},D_{24}),(A_{4m}\rtimes\ZZ_2,A_{4m-1},D_{8m})\}.$$
Next, we examine that $(X,G,D)$ can be $(A_{4m}\rtimes\ZZ_2,A_{4m-1},D_{8m})$, where $m\ge2$.
Let $$D=\lg a,b\di a=(1 2 \cdots 4m)(4m+1\,4m+2),~b=(1\,4m)(2\,4m-1)\cdots(2m\,2m+1)\rg.$$
Then we get $D\cong D_{8m}$, $b\in A_{4m}$ and $D\leq S_{4m+2}$.
Let $G=A_{4m-1}\leq S_m$ and $X=A_{4m}\rtimes\lg ab\rg$.
Then $X=A_{4m}\times\lg (4m+1\,4m+2)\rg\cong A_{4m}\times\ZZ_2$ and one yields $X=GD$.
Since $A_{4m}$ is simple for $m\ge2$, we know $D_X=1$.
Hence,
$$(X, G, D)=(\Aut(M_{12}),M_{11},D_{24})~{\rm or}~(A_{4m}\rtimes\ZZ_2,A_{4m-1},D_{8m}),$$
which implies that Rows 7 and 8 in Table~\ref{dihedral-skew} hold.

Suppose that $H$ is not quasiprimitive.
Then there exists a maximal subgroup $K$ of $H$ such that $K_X\neq1$ and $G\leq K$.
Clearly, $K$ is either $G\lg a_1\rg$ or $G\lg a_1,b_2\rg$, where $a_1\in\lg a^2\rg$ and $b_2\in\lg a^2,b\rg\setminus\lg a^2\rg$.
Suppose $K=G\lg a_1\rg$.
Then using Lemma~\ref{Ga} for $H$,  we get $K=G\lg a^2\rg\character H\lhd X$, which implies $K\lhd X$.
Hence, $G\lg a\rg\leq X$, contradicting to $\mathfrak{A}$.
Suppose $K=G\lg a_1,b_2\rg$.
Furthermore, we assume that $G\lg a^2\rg$ is not a subgroup of $H$.
Then with the same argument as Case 1, if $K_H\cap G=1$, then $G\lg a^2\rg\leq H$, which implies $G\lg a\rg\leq X$, a contradiction again.
Hence, $K_H\cap G=G$ and so $K_H$ is either $G\lg a_1\rg$ or $G\lg a_1^2,b_3\rg\leq K_H$, where $b_3\in\lg a_1,b_2\rg\setminus\lg a_1\rg$.
Suppose $K_H=G\lg a_1\rg$. Then $G\lg a^2\rg\leq H$, a contradiction.
Suppose $G\lg a_1^2,b_3\rg\leq K_H$.
Then $|K:K_H|<|H:K_H|=2$, which implies $K_H=K$.
Then we get $K=K_H=G\lg a^4,b_2\rg$ and $4\di \o(a)$. Clearly, $\lg a^4\rg\leq K_X$ and $b_2\notin K_X$.
Since $|X:K|=4$ and $X/K_X\leq S_4$, we know $G\leq K_X$, which implies $K_X=G\lg a^4\rg$.
Then we also get $G\lg a^2\rg\leq H$, a contradiction again.
\qed

\vskip 3mm
In the following lemma, we shall deal with the case that $D_X\neq1$.

\begin{lem}\label{GD1}
Suppose that $D_X\neq1$. Then $1<D_X<\lg a\rg$ and $X=D_X\rtimes Y$, where $(Y,G)$ is
  either $(PGL(2,11),A_5)$ or $(A_{2m+1}\rtimes\ZZ_2,A_{2m})$ with $m\ge 3$.
\end{lem}
\demo
By $D_X\unlhd D$, we get
$$D_X\in\{\lg a_1\rg, \lg a^2,b_1\rg, D\di 1\neq a_1\in\lg a\rg,b_1\in D\setminus\lg a\rg\}.$$
Since $X/C_X(D_X)\leq \Aut(D_X)$, a solvable group, and $G$ is a nonabelian simple group, we get $G\leq C_X(D_X)$. Hence, $\lg D_X, G\rg=G\times D_X$.
If $G\times D_X\unlhd X$, then $G\character (G\times D_X)\unlhd X$,
which implies $G\lhd X$, a contradiction.
Therefore, $G\times D_X\ntrianglelefteq X$, which implies $D_X=\lg a_1\rg<\lg a\rg$ with $a_1\neq1$,
and $X=(G\times\lg a_1\rg)D$.
Note that the quotient group $\ox:=X/\lg a_1\rg=\og\overline{D}$ is {\it PSDG}
where both $\og\cong G$ and $\overline{D}$ are core free and $\og\cap\overline{D}=1$.
Thus, by Lemma~\ref{GX1}, $\ox$ is one of Table~\ref{dihedral-skew}.
Set $\ox\cong Y$.
By $\lg a\rg\leq C_X(\lg a_1\rg)$, we get that $C_X(\lg a_1\rg)$ is either $X$ or $G\lg a\rg$.

Suppose that $C_X(\lg a_1\rg)=X$. Then by $a_1^{-1}=a_1^b=a_1\neq1$, we get $\o(a_1)=2$.
Then $X$ is a central extension of $\ZZ_2$ by $Y$, where $Y$ is one of groups in  Table~\ref{dihedral-skew}.
For a contradiction, assume that $X$ is a proper central extension of $\ZZ_2$ by $Y$.
Then by checking Magma, $\ox$ must not be $\AGL(3,2)$, $M_{12}$, $M_{24}$ or $\Aut(M_{12})$.
If $\ox=\PGL(2,11)$, then we get $\og\lg\ola\rg=\PSL(2,11)$ and so $G\lg a\rg$ is a proper central extension of $\ZZ_2$ by $\PSL(2,11)$. But it is impossible by Proposition~\ref{Schur}.(1), a contradiction.
With the same argument as the case $\PGL(2,11)$, we get $\ox\neq A_{2m_1+1}\rtimes\ZZ_2$ or $A_{4m_2}\rtimes\ZZ_2$ where $m_1\ge 3$ and $m_2\ge 2$.
Also,  by Proposition~\ref{Schur}.(1) again, we get $\ox\neq A_{4m}$ where $m\ge 2$.
Therefore, $X\cong\ZZ_2\times Y$ and so $\o(\ola)$ is odd.
Then $Y$ is either $\PSL(2,11)$ or $A_{2m+1}\rtimes\ZZ_2$, where $m\ge 3$, as desired.

Suppose that $C_X(\lg a_1\rg)=G\lg a\rg$.
Then $\ox=(\og\lg\ola\rg)\rtimes\lg \olb\rg$. Set $\og\lg\ola\rg\cong B$.
By Lemma~\ref{Ga}, we get
$(\ox,B,\og,\overline{D})=(\PGL(2,11), \PSL(2,11), A_5,D_{22})$ or
$(A_{m+1}\rtimes\ZZ_2, A_{m+1}, A_{m},D_{2(m+1)})$, with $m\ge6$ even.
Since the pre-image $\og\lg\ola\rg$ is a central extension of $\lg a_1\rg$ by $B$, by Proposition~\ref{Schur},
we get $G\lg a\rg=\lg a_1\rg\times(G\lg a_2\rg)$, where $\lg a\rg=\lg a_1\rg\times\lg a_2\rg$
and $G\lg a_2\rg$ is either $\PSL(2,11)$ or $A_{m+1}$.
By $G\lg a_2\rg\character G\lg a\rg\lhd X$, we know $G\lg a_2\rg\lhd X$.
Then $G\lg a_2,b\rg$ is the complement group of $\lg a_1\rg$ in $X$.
Then $$X=\lg a_1\rg\rtimes G\lg a_2,b\rg,$$
where $G\lg a_2,b\rg$ is either $\PGL(2,11)$ or $A_{m+1}\rtimes\ZZ_2$ with $m\ge6$ even, as desired.
\qed

\subsection{$G\cap D=\ZZ_2$}\label{subsec:GD2}
In this subsection, let $X=GD$ be {\it PSDG} such that $G\cap D=\ZZ_2$ and $G_X=1$.
By Lemma~\ref{D-goup}, we get $G\cap D=\lg b_1\rg$, with $b_1\in D\setminus\lg a\rg$.
We may as well choose notations so that $b_1 = b$.
Then $X=G\lg a\rg$, where $G\cap\lg a\rg=1$, $G$ is core free and there exits an involution $b\in G$ such that $a^b=a^{-1}$.

\begin{lem}\label{GC1}
Theorem~\ref{main1}.(2) holds, that is, $(X, G, D)$ is either $(A_{m+1}, A_m,D_{2(m+1)})$ or $(A_{m+1}\times\ZZ_2, A_m,D_{4(m+1)})$, with $\frac m2 \ge 3$ even.
\end{lem}
\demo
Consider $\lg a\rg_X$, the core of $\lg a\rg$ in $X$.
Suppose that $\lg a\rg_X=1$. Then by Proposition~\ref{skew}, we get
$$(X, G, \lg a\rg)\in\{(\PSL(2,11), A_5,\ZZ_{11}),(M_{23}, M_{22},\ZZ_{23}),(A_{m+1}, A_m,\ZZ_{m+1})\},$$ where $m \ge 6$ is even.
Note that $D=\lg a,b\rg\leq X$ where $b\in G$.
Using Magma~\cite{BCP1997}, we know that there no exists a subgroup isomorphic to $D_{22}$ (resp. $D_{46}$) in $\PSL(2,11)$ (resp. $M_{23}$), which implies that $(X, G, \lg a\rg)\notin\{(\PSL(2,11), A_5,\ZZ_{11}),(M_{23}, M_{22},\ZZ_{23})\}$.
Next, suppose that $(X,G,\lg a\rg)=(A_{m+1}, A_m,\ZZ_{m+1})$, where $m\ge6$ is even.
Then we can choose $a=(1\,2\,\cdots\,m\,m+1).$
Clearly, there exists an involution $b=(1\,m+1)(2\,m)\cdots(\frac m2\,\frac m2+2)$ in $S_{m+1}$ such that
$\lg a, b\rg\cong D_{2(m+1)}$.
Moreover, $b\in A_{m+1}$ if and only if $\frac m2$ is even.
Then we get
$$(X, G, D)=(A_{m+1}, A_m,D_{2(m+1)}),~ {\rm with}~ 4\di m~{\rm and}~m\ge 8.$$
So in what follows, we assume that $\lg a\rg_X\neq1$.

Set $\lg a\rg_X=\lg a_1\rg$ with $1\neq a_1\in\lg a\rg$.
Since $X/C_X(\lg a_1\rg)\leq\Aut(\lg a_1\rg)$ is abelian and $\lg a\rg\leq C_X(\lg a_1\rg)$,
we get $C_X(\lg a_1\rg)=G\lg a\rg=X$, which implies $b\in C_X(\lg a_1\rg)$, and so $\o(a_1)=2$.
If $G\times\lg a_1\rg\lhd X$, then by $G\character(G\times\lg a_1\rg)\lhd X$, we get $G\lhd X$, a contradiction.
Consider $\ox=X/\lg a_1\rg=\og\lg\ola\rg$.
Since $\og\cap\lg\ola\rg=1$ and $\og_{\ox}=\lg\ola\rg_{\ox}=1$,
by the result as the above, we know $\ox=A_{m+1}$, with $4\di m$ and $m\ge 8$.
Hence, $X$ is a central extension of $\lg a_1\rg\cong\ZZ_2$ by $A_{m+1}$, with $4\di m$ and $m\ge 8$.
Then $(X, G, D)=(\ZZ_2.A_{m+1}, A_m,D_{2(m+1)})$ with $4\di m$ and $m\ge 8$.
By Proposition~\ref{Schur}, we know
$$(X, G, D)=(A_{m+1}\times\ZZ_2, A_m,D_{4(m+1)}),~ {\rm with}~ 4\di m~{\rm and}~m\ge 8,$$
as desired.
\qed
\section{Further  research}

Since all nonabelian simple groups are finite nonabelian characteristically simple groups, we propose the following problem for further researches.
\begin{que}
Determine the product of finite nonabelian characteristically simple groups and dihedral groups.
\end{que}

\vskip 3mm
\f {\bf Acknowledgments:}
The author would like to thank the reviewers for their valuable comments and suggestions.

\section{Appendices}
Magma Code for [Corollary 2.2]

\begin{lstlisting}[language=Magma, escapechar=@]
G:=Alt(6);
L:=CyclicGroup(2);
L1:=DihedralGroup(8);
X:=SolvableQuotient(L1);
E := ExtensionsOfElementaryAbelianGroup(2,4,G);
for i:=1 to 1 do
E1:=E[i];
_,E4:=CosetAction(E1,sub<E1| Id(E1)>);
E2:=DirectProduct(E4,L);
_,E3:=CosetAction(E2,sub<E2| Id(E2)>);
A:=Subgroups(E3:OrderEqual:=16);
for j:=1 to #A do
B:=A[j]@\textasciigrave@subgroup;
Y:=Core(E3,B);
if (#Y eq 1) then
Y:=SolvableQuotient(B);
if (IsIsomorphic(Y,X)) then
[i,j];
end if;
end if;
end for;
end for;
\end{lstlisting}

\vskip 3mm
\begin{lstlisting}[language=Magma, escapechar=@]
G:=Alt(7);
L:=CyclicGroup(2);
L1:=DihedralGroup(8);
X:=SolvableQuotient(L1);
E := ExtensionsOfElementaryAbelianGroup(2,4,G);
for i:=1 to #E do
E1:=E[i];
_,E4:=CosetAction(E1,sub<E1| Id(E1)>);
E2:=DirectProduct(E4,L);
_,E3:=CosetAction(E2,sub<E2| Id(E2)>);
A:=Subgroups(E3:OrderEqual:=16);
for j:=1 to #A do
B:=A[j]@\textasciigrave@subgroup;
Y:=Core(E3,B);
if (#Y eq 1) then
Y:=SolvableQuotient(B);
if (IsIsomorphic(Y,X)) then
[i,j];
end if;
end if;
end for;end for;
\end{lstlisting}

\vskip 3mm
\begin{lstlisting}[language=Magma, escapechar=@]
G:=AGL(4,2);
K1:=GL(4,2);
S:=Subgroups(G:OrderEqual:=#K1);
[[i,S[i]@\textasciigrave@order] : i in [1..#S]];
H2:=S[1]@\textasciigrave@subgroup;

H:=DihedralGroup(8);
X:=SolvableQuotient(H);
S:=Subgroups(G:OrderEqual:=16);
#S;
for i:=1 to #S do
L1:=S[i]@\textasciigrave@subgroup;
Y:=SolvableQuotient(L1);
if (IsIsomorphic(Y,X)) then
i;
Core(G,L1);
H2 meet L1;
end if;
end for;
\end{lstlisting}

\vskip 3mm
\begin{lstlisting}[language=Magma, escapechar=@]
G1:=AGL(3,2);
_,G:=CosetAction(G1,sub<G1| Id(G1)>);
K2:=GL(3,2);
_,K:=CosetAction(K2,sub<K2| Id(K2)>);
H:=DihedralGroup(4);
X:=SolvableQuotient(H);
S:=Subgroups(G:OrderEqual:=#K);
A:=Subgroups(G:OrderEqual:=#X);
for j:=1 to #A do
L1:=A[j]@\textasciigrave@subgroup;
Y:=SolvableQuotient(L1);
if (IsIsomorphic(Y,X)) then
Z:=Core(G,L1);
if (#Z eq 1) then
for i:=1 to #S do
H2:=S[i]@\textasciigrave@subgroup;
_,Y:=CosetAction(H2,sub<H2| Id(H2)>);
B:=H2 meet L1;
if (IsIsomorphic(Y,K) and #B eq 1) then
[i,j];
end if;
end for;
end if;
end if;
end for;
\end{lstlisting}

\vskip 3mm
\begin{lstlisting}[language=Magma, escapechar=@]
load "M12";
#G;
S:=Subgroups(G:OrderEqual:=12);
[[i,S[i]@\textasciigrave@order] : i in [1..#S]];
H:=DihedralGroup(6);
X:=SolvableQuotient(H);
for i:=1 to #S do
H1:=S[i]@\textasciigrave@subgroup;
Y:=SolvableQuotient(H1);
if (IsIsomorphic(Y,X)) then
i;
Core(G,H1);
end if;
end for;

T:=Subgroups(G:OrderEqual:=7920);
[[i,T[i]@\textasciigrave@order] : i in [1..#T]];
K1:=T[1]@\textasciigrave@subgroup;

for i:=6 to 8 do
H1:=S[i]@\textasciigrave@subgroup;
K1 meet H1;
end for;
\end{lstlisting}

\vskip 3mm
\begin{lstlisting}[language=Magma, escapechar=@]
load "M24";
#G;
S:=Subgroups(G:OrderEqual:=48);
[[i,S[i]@\textasciigrave@order] : i in [1..#S]];
H:=DihedralGroup(24);
X:=SolvableQuotient(H);
for i:=1 to #S do
H1:=S[i]@\textasciigrave@subgroup;
Y:=SolvableQuotient(H1);
i;
IsIsomorphic(Y,X);
end for;
\end{lstlisting}

\vskip 3mm

Magma Code for [Lemma 2.5]

\begin{lstlisting}[language=Magma, escapechar=@]
G:=PGL(2,11);
G;
H:=DihedralGroup(11);
X:=SolvableQuotient(H);
S:=Subgroups(G:OrderEqual:=60);
[[i,S[i]@\textasciigrave@order] : i in [1..#S]];
H1:=S[1]@\textasciigrave@subgroup;
H1;

S:=Subgroups(G:OrderEqual:=22);
[[i,S[i]@\textasciigrave@order] : i in [1..#S]];
L1:=S[1]@\textasciigrave@subgroup;
L1;
Y:=SolvableQuotient(L1);
IsIsomorphic(Y,X);
Core(G,L1);
H1 meet L1;
\end{lstlisting}

\vskip 3mm
\begin{lstlisting}[language=Magma, escapechar=@]
G:=PSL(2,11);
G;
L:=CyclicGroup(2);
K:=DirectProduct(G,L);
H:=DihedralGroup(11);
X:=SolvableQuotient(H);
S:=Subgroups(K:OrderEqual:=22);
[[i,S[i]@\textasciigrave@order] : i in [1..#S]];
for i:=1 to #S do
H1:=S[i]@\textasciigrave@subgroup;
Y:=SolvableQuotient(H1);
i;
IsIsomorphic(Y,X);
end for;
\end{lstlisting}

\vskip 3mm
\begin{lstlisting}[language=Magma, escapechar=@]
load "M12";
#G;
L:=CyclicGroup(2);
K:=DirectProduct(G,L);
H:=DihedralGroup(12);
X:=SolvableQuotient(H);
S:=Subgroups(K:OrderEqual:=24);
for i:=1 to #S do
H1:=S[i]@\textasciigrave@subgroup;
Y:=SolvableQuotient(H1);
if (IsIsomorphic(Y,X)) then
i;
end if;
end for;
\end{lstlisting}

\vskip 3mm
\begin{lstlisting}[language=Magma, escapechar=@]
load "M12";
#G;
S:=Subgroups(G:OrderEqual:=24);
[[i,S[i]@\textasciigrave@order] : i in [1..#S]];
H:=DihedralGroup(12);
X:=SolvableQuotient(H);
for i:=1 to #S do
H1:=S[i]@\textasciigrave@subgroup;
Y:=SolvableQuotient(H1);
i;
IsIsomorphic(Y,X);
end for;
K:=PermutationGroup(AutomorphismGroup(G));
#K;
S:=Subgroups(K:OrderEqual:=24);
[[i,S[i]@\textasciigrave@order] : i in [1..#S]];
for i:=1 to #S do
H1:=S[i]@\textasciigrave@subgroup;
Y:=SolvableQuotient(H1);
i;
IsIsomorphic(Y,X);
end for;
H1:=S[12]@\textasciigrave@subgroup;
Core(K,H1);
H2:=S[18]@\textasciigrave@subgroup;

S:=Subgroups(K:OrderEqual:=7920);
[[i,S[i]@\textasciigrave@order] : i in [1..#S]];
K1:=S[1]@\textasciigrave@subgroup;

K1 meet H1;
K1 meet H2;
\end{lstlisting}

\vskip 3mm
\begin{lstlisting}[language=Magma, escapechar=@]
G2:=AGL(3,2);
L:=CyclicGroup(2);
G1:=DirectProduct(G2,L);
_,G:=CosetAction(G1,sub<G1| Id(G1)>);
H:=DihedralGroup(8);
X:=SolvableQuotient(H);
S:=Subgroups(G:OrderEqual:=16);
for i:=1 to #S do
L1:=S[i]@\textasciigrave@subgroup;
Y:=SolvableQuotient(L1);
if (IsIsomorphic(Y,X)) then
i;
Core(G,L1);
end if;
end for;
\end{lstlisting}

\vskip 3mm
\begin{lstlisting}[language=Magma, escapechar=@]
load "M24";
#G;
L:=CyclicGroup(2);
K:=DirectProduct(G,L);
S:=Subgroups(K:OrderEqual:=48);
H:=DihedralGroup(24);
X:=SolvableQuotient(H);
for i:=1 to #S do
H1:=S[i]@\textasciigrave@subgroup;
Y:=SolvableQuotient(H1);
if (IsIsomorphic(Y,X)) then
i;
end if;
end for;
\end{lstlisting}

\end{document}